\documentclass[12pt]{amsart}
\usepackage{amsmath,amsthm,epsfig,amssymb}
\usepackage{amsthm,amsfonts,amssymb,amscd}
\usepackage[all]{xy}
\newtheorem{thm}[subsection]{Theorem}
\newtheorem{prop}[subsection]{Proposition}

\newtheorem{lem}[subsection]{Lemma}

\newtheorem{rem}[subsection]{Remark}
\theoremstyle{definition}
\newtheorem{Def}[subsection]{Definition}

\newtheorem{exam}[subsection]{Example}
\newtheorem{proposition-definition}[subsection]{Proposition-Definition}

\newcommand{\CC}{{\mathbb C}}

\newcommand{\ZZ}{{\mathbb Z}}

\newcommand{\PP}{{\mathbb P}}

\newcommand{\OOO}{{\mathcal O}} 
\newcommand{\III}{{\mathcal I}}

\newcommand{\CCC}{{\mathcal C}}

\author{A. El Mazouni}

\address{Laboratoire de Math\'ematiques de Lens EA 2462
Facult\'e des Sciences Jean Perrin
Rue Jean Souvraz, SP18
F-62307 LENS  Cedex France}
\email{mazouni@euler.univ-artois.fr}

\author{F. Laytimi}

\address{F. L.: Math\'ematiques - b\^{a}t. M2, Universit\'e Lille 1,
F-59655 Villeneuve d'Ascq Cedex, France}
\email{fatima.laytimi@math.univ-lille1.fr}

\author{D.S. Nagaraj}

\address{Institute of Mathematical Sciences C.I.T. campus, Taramani, 
chennai 600113,India}
\email{dsn@imsc.res.in}

\title{Morphisms from $\PP^2$ to $Gr(2, \CC^{4}).$}
\date{}
\begin{document}

\begin{abstract}  
In this note we study morphisms from $\PP^2 $ to $ Gr(2, \CC^{4})$ from
the point of view of the cohomology class they represent in the Grassmannian.
This leads to some new result about projection of $d$-uple imbedding of 
$\PP^2$  to $\PP^5.$
\end{abstract}

\maketitle

\section{Introduction} \setcounter{page}{1}
We denote by $\PP^2$  the projective plane over the field of complex numbers and 
$Gr(2, \CC^{4})$ be the Grassmannian of two-dimentional quotient spaces of $\CC^{4}.$
In this paper we investigate the  possible types  of non-constant morphisms 
$\PP^2  \to  Gr(2, \CC^{4}).$
Any non-constant morphism $\PP^2  \to  Gr(2, \CC^{4})$ determines a 
cohomology class in $H^4(Gr(2, \CC^{4}), \ZZ).$ We consider the 
following problem:\\
Determine the necessary and sufficient conditions that a  cohomology class of $H^4(Gr(2, \CC^{4}), \ZZ)$ 
has to satisfy in order to be  represented by a morphism  from $\PP^2 \to  Gr(2, \CC^{4})$? 
 
It is easy to see that if a cohomology class is  represented by a morphism from $\PP^2$  then it has  to 
satisfy an obvious  necessary condition (see Lemma (\ref{le21})).
We show that in general this condition is not sufficient.
The following result shows that  there are classes satisfying the condition but are 
not represented by  morphisms:

\begin{thm} Let $c$ and $a$ be two integers.
  Assume  $c \geq 4 $  and $1\leq a \leq c-2$ or $c^2-c+2 \leq a \leq c^2-1.$
 Then the cohomology class $(a,c^2-a)$  is not represented by a morphism 
 $ \PP^2 \to Gr(2, \CC^{4}).$  
 
\end{thm}

 The following theorem shows that in some cases  morphisms exist:
 
 \begin{thm}
 0) For every $c\geq 1$ the cohomology classes $(0,c^2)$ and $(c^2,0)$
 are represented by  morphisms $\PP^2 \to Gr(2, \CC^{4}).$
  
 1) Let $1\leq c \leq 3$ and $0\leq \ell \leq c^2$ be integers.
 Then  the cohomology class $(\ell,c^2-\ell)$
 is represented by a morphism  $\PP^2 \to Gr(2, \CC^{4}).$
 
 2)  If $c=4$ and $3\leq \ell \leq 13$ then there are morphisms 
$\PP^2 \to Gr(2, \CC^{4})$ representing the cohomology class $(\ell, 16-\ell).$

 3) Let $c \geq 5 $ be an integer. Let $k$ be the largest integer such that 
 $k.c \leq (c^2-3c+2)/2.$ Then for every integer $\ell$ in one of the following
 intervals  $[t(c-3)+2, t.c]$ for $1\leq t \leq k$ or $[(c^2-3c+2)/2+1, c^2/2]$
 there are morphisms $ \PP^2 \to Gr(2, \CC^{4})$ representing the  
 cohomology class $(\ell,c^2-\ell).$ 
 Also, for every such $\ell$ there is a morphism whose cohomology class 
 is $(c^2-\ell, \ell).$
  \end{thm}

It is also shown here (see, Remark (\ref{re23})) that for every integer $n \geq 1$ 
there are morphisms  
$$f_n: \PP^2  \to  Gr(2, \CC^{4})$$
 such that $f_n$ is one to one onto its image and $f_n^*(\OOO_{Gr(2,\CC^4)}(1)) = 
\OOO_{\PP^2}(n).$ 

{\em Acknowledgments:}  The third author would like to thank, IFIM,  
Universit\'e Lille,  and Universit\'e d' Artois, France.  

\section{Morphisms to Grassmannian}
Here we recall some  results about morphisms  from a variety to a Grassmannian 
variety.

\begin{Def} Let $X$ be a variety over the field $\CC$ of complex numbers and
$s$ be a positive integer. 
A vector bundle $E$ over $X$ is said to be generated by $s$ sections if there a
surjection of vector bundles
$$\CC^{s} \otimes \OOO_X \rightarrow  E.$$ 
\end{Def}

\begin{Def} Let $X$ be a projective variety. Let $r$ and $k$ be
two positive integers. Let $E_1$ and $E_2$ be two vector bundles
of rank $r$ over $X.$  Two vector bundle  surjections 
$$\phi_1: \CC^{r+k} \otimes \OOO_X \rightarrow  E_1,$$
and
$$\phi_2: \CC^{r+k} \otimes \OOO_X \rightarrow  E_2,$$
are said to be equivalent if there exists an isomorphsim  of 
vector bundles  $\psi: E_1 \to E_2$ over $X$ such that the 
following diagram commutes: 
$$
\begin{array}{ccl}
  \CC^{r+k}\otimes \OOO _X &\stackrel{\phi_1}{\longrightarrow} & E_1\\
  \Vert & {} & \downarrow{\psi} \\ 
  \CC^{r+k}\otimes \OOO_X &\stackrel{\phi_2}{\longrightarrow} & E_2 
\end{array} 
$$ 
\end{Def}
 
The following two lemmas are well known(see, for example \cite{FD}).
We include the proofs of these lemmas for the sake of completeness.

\begin{lem}\label{cors} Let $X$ be a projective variety. Let $r$ and $k $ be
two positive integers. There is a natural bijection between the following two
sets:

1) The set of equivalence classes of surjection of vector bundles
$$\phi: \CC^{r+k} \otimes \OOO_X \rightarrow  E,$$
where $E$ is a vector bundle of rank $r$ on $X.$

2) The set of morphisms $f: X \rightarrow Gr(r,\CC^{r+k})$, where
$Gr(r,\CC^{r+k})$ is the Grassmannian of r dimensional quotient of 
$\CC^{r+k} .$
\end{lem}

{\it Proof:}  
Given a surjection $\phi: \CC^{r+k} \otimes \OOO _X \rightarrow E$ by
sending $$x\mapsto \{\CC^{r+k} \rightarrow E_x\}$$ defines a
morphism $f:X \longrightarrow Gr(r, \CC^{r+k})$.  This defines a map
from the set in 1) to the set in 2).

To prove the existence of the map in the other direction, we first note the
following: on 
$Gr(r, \CC^{r+k})$ there is a canonical vector bundle surjection 
$\CC^{r+k} \otimes \OOO_{Gr(r, \CC^{r+k})} \longrightarrow Q,$ where
$Q$ is the rank $r$ bundle on $Gr(r, \CC^{r+k})$ whose fiber at $x \in
Gr(r, \CC^{r+k})$ is the quotient vector space corresponding to $x.$

Now given a morphism $f:X \longrightarrow Gr(r, \CC^{r+k})$ by pulling
back by $f$ the canonical surjection of vector bundles on the
Grassmannian, we get a surjection of vector bundles
$$
\CC^{r+k} \otimes \OOO _X =
f^*(\CC^{r+k} \otimes \OOO _{Gr(r, \CC^{r+k})})\longrightarrow f^*(Q).
$$ 

This gives a map from the set in 2) to the set in 1). 
These two maps are clearly inverse to each other hence we get the 
required bijection.  $\hfill{\Box}$

\begin{thm}\label{thm1} Let $X$ be a projective variety.  There is a natural
  bijection between the following two sets:

1) The set of equivalence classes of surjections of vector bundles
$$\CC^{r+k} \otimes \OOO _X\longrightarrow E$$
with $rank(E)=r$ and $\det(E)$ is ample.

2) The set of morphisms $f:X \longrightarrow Gr(r, \CC^{r+k})$ with
$f$ finite (onto its image).
\end{thm}

{\it Proof:} Given an element of the set in 1), i.e., a surjection, $\CC^{r+k}
\otimes \OOO _X\longrightarrow E$ with $\rm{rk}(E)=r$ and $\det(E)$ is
ample, then as in Lemma(\ref{cors}) it determines a morphism $f:X
\longrightarrow Gr(r,\CC^{r+k})$, also, it follows from Lemma(\ref{cors}), 
that $f^*(Q)\simeq E$.  Hence $f^*{\rm det}(Q)\simeq \rm{det} (E)$. 
Since
${\rm det}(E)|_{f^{-1}(f(x))} \simeq f^*\rm{det}(Q) |_{f^{-1}(f(x))}$ 
is trivial, the ampleness assumption on ${\rm det}(E)$ 
implies $\rm{dim}({f^{-1}(f(x))})\leq 0.$ Thus, if $\det(E)$ is ample,
then $f$ is finite onto its image. 

In the other direction, if $f:X \longrightarrow Gr(r,\CC^{r+k})$ is
finite morphism onto its image, then
$$\rm{det} (f^*Q)\simeq f^*(\rm{det} Q)$$ is ample, as pull-back
of an ample bundle remains ample under finite morphism. 
 $\hfill{\Box}$ 

Let $L$ be an vector subspace of dimension one in $\CC^4.$ Then the projective
$2$ plane 
$$ \PP_L^2 = Gr(2, {\CC^4}/L) $$ 
is naturally a subvariety  of $Gr(2, \CC^{4}).$  Similarly, if $H$ is an one dimensional quotient vector
space of $\CC^4$, then the  projective $2$ plane 
$$\PP_H^2 = Gr(2, {(\CC^4)^*}/{H^*}) $$ 
is also naturally a subvariety of $Gr(2, \CC^{4}),$ where $(\CC^4)^*$ (resp. $H^*$) 
denotes the dual vector space of $\CC^4 $ (resp. of $H$).
Note that $\PP_H^2$ is the subset of $Gr(2, {\CC^4})$ consists of all those
quotient vector spaces of dimension two of $\CC^{4}$ which has $H$ as their quotient.

\begin{Def}\label{triv}  A morphism $f: \PP^2 \to Gr(2, \CC^{4})$ is said to be trivial
if it is one of the following:
\begin{itemize}
\item a) $f$ is a constant  morphism.
\item b) Image of $f$ is $\PP_L^2$ for some one dimensional subspace $L $ of $ \CC^4.$ 
\item c) Image of $f$ is $\PP_H^2$ for some one dimensional quotient $ H$ of $\CC^4.$
\end{itemize}
\end{Def}

\begin{exam}\label{ex1}
It is known that (see, for example \cite{SD}) the Veronese  surface (2-uple embedding of $\PP^2$ in 
$\PP^5$) is contained in a smooth quadric in $\PP^5.$  As any two smooth quadrics
are isomorphic via a projective automorphism we see that there is an embedding
of $\PP^2$ in $Gr(2, \CC^{4}).$ In fact, the global sections
$$ s_1=(X,Z),\ \  s_2=(Y,X), \ \ s_3=(Z,Y) \  \, {\rm and} \ \, s_4=(X,0) $$
of $ \OOO_{\PP^2}(1) \oplus \OOO_{\PP^2} (1)$ gives a  vector bundle surjection
$$\CC^{4} \otimes \OOO_{\PP^2} \rightarrow  \OOO_{\PP^2}(1) \oplus \OOO_{\PP^2} (1).$$
Then by Lemma (\ref{cors}) we get morphism $\phi : \PP^2 \to Gr(2, \CC^{4}).$ It is
easy to see that $\phi$ is an embedding and $\phi$ composed with the natural embedding
of $Gr(2, \CC^{4})$ into $\PP^5$ is  given by the quadrics 
$$X^2-YZ, \ \ XY-Z^2, \ \ Y^2-XZ, \ \ XZ, \ \ X^2, \ \ XY.$$
As these quadrics form a basis for the space of quadrics on $\PP^2$ we get that
the veronese  embedding of $\PP^2$ in $\PP^5$ factors through $Gr(2, \CC^{4}).$ 
\end{exam}

\begin{exam}\label{ex2}
 More generally, for any positive integers $a$ and $b$ 
the global sections
$$ s_1=(X^a,Z^b), \ \ s_2=(Y^a,X^b),\ \ s_3=(Z^a,Y^b) \ \, {\rm and} \ \, s_4=(X^a,0) $$
of $ \OOO_{\PP^2}(a) \oplus \OOO_{\PP^2} (b)$ gives a  vector bundle surjection
$$\CC^{4} \otimes \OOO_{\PP^2} \rightarrow  \OOO_{\PP^2}(a) \oplus \OOO_{\PP^2} (b).$$
Then by Lemma (\ref{cors}) we get a morphism $\phi : \PP^2 \to Gr(2, \CC^{4}).$
\end{exam}

\begin{exam} \label{ex3}
Let $T_{\PP^2}$ denotes the tangent bundle of $\PP^2.$ Then one has an exact
sequence (See, page 409 of \cite{GH})
$$ 0 \to \OOO_{\PP^2} \to \OOO_{\PP^2}(1)^3 \to  T_{\PP^2} \to 0.$$
>From this we see that the space of all global sections ${\rm H}^0(T_{\PP^2})$
of the tangent bundle is a vector space of dimension eight.  Now it is easy to see
that if we choose four linearly independent general sections of the tangent
bundle we get a surjective morphism 
$$ \OOO_{\PP^2}^4 \to  T_{\PP^2}. $$
This surjection give rise to a morphism from $\PP^2$  to $Gr(2, \CC^{4}).$
Since the tangent bundle is not direct sum of line bundles, this morphism is 
different from the morphisms given by example above.
\end{exam}
  
\section{Non-trivial morphisms from $\PP^2$  to $Gr(2, \CC^{4})$}
 
 On the Grassmaian  $Gr(2, \CC^{4})$ one has universal exact sequence:
\begin{equation}\label{univ2}
 0 \to S \to \CC^{4} \otimes \OOO_{Gr(2, \CC^{4})} \to Q \to 0, 
\end{equation}
 where $S$ and $Q$ are respectively the universal sub bundle and quotient bundle of 
rank two on $Gr(2, \CC^{4}).$ The fiber of $Q$ (resp. of $S$) at a point $p \in Gr(2, \CC^{4})$
is the two dimensional quotient space (resp. subspace, which is the kernal of this quotient map)
of $\CC^{4}$ corresponding to the point $p.$
It is known (see, page 197 and 411 of \cite{GH}) that the cohomology group $H^4(Gr(2, \CC^{4}), \ZZ) $
is equal to 
\begin{equation}\label{c2}
 \ZZ[c_2(Q)] \oplus  \ZZ[c_2(S)], 
 \end{equation}
 where $c_2(Q)$ (resp. $c_2(S)$) is the second Chern class of $Q$ (resp. of $S$).
 If $ \phi : \PP^2 \to Gr(2, \CC^{4})$ is a non constant morphism then the cohomology class of 
 $\phi_{*}([\PP^2])$ is an element of  $H^4(Gr(2, \CC^{4}), \ZZ). $
 It is easy to see that the morphism b)  (resp. c)) of  Example(\ref{triv}) gives the cohomology class $(0,1)$ 
 (resp. $(1,0)$) of the decomposition in \eqref{c2} of the cohomology group. The exact sequence 
 corresponding to cohomology class $(0,1)$ is 
 \begin{equation}\label{trivialclass}
 0 \to \Omega^1(1) \to \CC^{4} \otimes \OOO_{\PP^2} \to \OOO_{\PP^2}(1) \oplus \OOO_{\PP^2} \to 0. 
\end{equation}
The dual of the exact sequence \eqref{trivialclass} correspondence to the cohomology class $(1,0).$  
 
 {\bf Question 1)}  a) Given the cohomology class $(a,b)$ of the decomposition in
  \eqref{c2} of the cohomology group does there exists a morphism 
  $\phi :  \PP^2 \to Gr(2, \CC^{4})$ such that the cohomology class of this morphism is $(a,b)?$
  
  b) For which cohomology class $(a,b)$ does there exists a  generically injective morphism 
  $\phi :  \PP^2 \to Gr(2, \CC^{4})$ such the cohomology class of this morphism is $(a,b)?$ 

{\bf Question 2)}  Let $Q$ be the vector bundle on $Gr(2, \CC^{4})$ as in the equation 
\eqref{univ2}. For  which cohomology classes $(a,b)$   does there exists a 
morphism   $\phi :  \PP^2 \to Gr(2, \CC^{4})$  the bundle $\phi^*(Q)$ is indecomposable?

\begin{rem}\label{re-1} Let $\phi :  \PP^2 \to Gr(2, \CC^{4})$ be a morphism.
 By pulling back the universal exact sequence
\eqref{univ2} on $Gr(2, \CC^{4})$ we get an exact sequence 
\begin{equation}\label{univ3}
 0 \to \phi^*(S) \to \CC^{4} \otimes \OOO_{\PP^2} \to \phi^*(Q )\to 0, 
\end{equation}
of vector bundles on $\PP^2.$ By dualizing the exact sequence \eqref{univ3} we get
another exact sequence 
\begin{equation}\label{univ6}
 0 \to \phi^*(Q)^{\vee} \to \CC^{4} \otimes \OOO_{\PP^2} \to \phi^*(S )^{\vee}\to 0,  
\end{equation}
of vector bundles on $\PP^2.$ By Lemma (\ref{cors}) the surjection 
\begin{equation}
\CC^{4} \otimes \OOO_{\PP^2} \to \phi^*(S )^{\vee}\to 0,  
\end{equation}
gives a morphism from $\PP^2 $ to $Gr(2, \CC^{4})$ which we denote by
$\phi^{\vee}$ and call the dual morphism. If  $(a,b)$
is the cohomology class of a morphism $\phi$ then it is clear that the cohomology class of 
the dual morphism  $\phi^{\vee}$ is $(b,a).$
\end{rem}

\begin{rem}\label{re0} 
Let $\phi : \PP^2  \to Gr(2, \CC^{4})$ be a
closed immersion.  Assume that $\phi$ is 
a non trivial imbedding  as in  Definition(\ref{triv}).   
It follows from  \cite{Ta}  the morphism   
 $\phi$ or $\phi^{\vee},$ is given by a surjection
 \begin{equation}\label{surj} 
\CC^{4} \otimes \OOO_{\PP^{2}} \to \OOO_{\PP^{2}}(1)\oplus \OOO_{\PP^{2}}(1).
\end{equation}
If  the morphism $\phi$ is given by  \eqref{surj} then
$\phi_*([\PP^2]) = (1,3)$ and  the class of  
$(\phi^{\vee})_*([\PP^2]) = (3,1).$ Thus the only classes $(a,b)$ represented by regularly  imbedded 
$\PP^2$ in $Gr(2, \CC^{4})$ are $(1,0), (0,1), (1,3)$ and $(3,1).$
\end{rem}

\begin{lem}\label{le21}  Let  
$\phi :  \PP^2 \to Gr(2, \CC^{4})$ 
be a non constant morphism. 
Then the cohomolgy class given by the morphism 
$\phi$ is of the form $(a, c^2-a),$ for some integers $c > 0$ and $0\leq a \leq {c^2}.$ 
\end{lem}

{\it Proof:} By pulling back the universal exact sequence \eqref{univ2} by $\phi$
we get the following  exact sequence of vector bundles on $\PP^2:$
\begin{equation}\label{exa}
 0 \to \phi^*(S) \to \CC^{4} \otimes \OOO_{\PP^{2}} \to \phi^*(Q) \to 0. 
\end{equation}
 Since, the bundle  
$\phi^*(Q)$ is non-trivial and is a quotient of trivial bundle,  we see that the Chern classes $c_1(\phi^*(Q))$
and $c_2(\phi^*(Q))$ are both non-negative integers (we identify 
$H^4(\PP^{2}, \ZZ)$ with $\ZZ$ by sending the class $[H].[H]$ to $1$ in $\ZZ,$ where $[H]$ is
the class of a hyper plane).   Set $a = c_2(\phi^*(Q))$ and $c = c_1(\phi^*(Q)).$
Note that  $c^2-a=c_2(\phi^*(S))= c_2((\phi^*(S))^{\vee}),$ where
$(\phi^*(S))^{\vee}$ is the dual of the bundle $\phi^*(S).$
Since $(\phi^*(S))^{\vee}$ is generated by sections we must have $c^2-a\geq 0.$ This completes the 
proof of the lemma.  $\hfill{\Box}$

\begin{rem}\label{re} The following question arises naturally.
Given $(a, c^2-a),$ with $c > 0$ and $0\leq a \leq {c^2}$
does there exists a morphism $\phi :  \PP^2 \to Gr(2, \CC^{4})$ such that 
the cohomolgy class given by the morphism 
$\phi$ is of the form $(a, c^2-a)?$  In general this question has negative answer 
(see Lemma (\ref{re1})
below). Also, we give below partial  answer to the above question 
(see Theorem(\ref{nonexist}) and 
Theorem(\ref{exist})). We believe that there are no morphisms from $\PP^2$
which represent the remaining cohomology classes.
\end{rem}
We need the following theorem:

\begin{thm}\label{CB}({\rm Cayley-Bacharach theorem})
Let $S$ be a non-singular surface and $Z$ be a sub scheme of S of dimension zero.
Then there exists a rank two vector bundle $E$ on $S$ with a section $s$ such
that the zero sub scheme of $s$ is $Z$ if and only if for every point $p\in Z$ the 
linear system 
$$|\III_ZK_S\otimes{\rm det}(E)| = |\III_{Z-\{p\}}K_S\otimes{\rm det}(E)|,$$  
where
$\III_Z$ (respectively, $\III_{Z-\{p\}}$) denotes the ideal sheaf of $Z$ (respectively, 
of  $Z-\{p\}).$
\end{thm}

{\it Proof:} This is a consequence of [Theorem 7.\cite{TV}].
(See, \cite{GH} page 731, for the case $Z$ is reduced).

\begin{lem}\label{re1} There is no morphism $\phi :  \PP^2 \to Gr(2, \CC^{4})$ such that 
the cohomolgy class given by the morphism 
$\phi$ is of the form $(1, 15).$ 
\end{lem}

{\it Proof:} Assume there exists a morphism $\phi :  \PP^2 \to Gr(2, \CC^{4})$ such that 
the cohomolgy class given by the morphism 
$\phi$ is of the form $(1, 15).$
By pulling back the universal exact sequence 
\eqref{univ2} by $\phi$
gives the following  exact sequence of vector bundles on $\PP^2:$
\begin{equation}\label{univ4}
 0 \to \phi^*(S) \to \CC^{4} \otimes \OOO_{\PP^{2}} \to \phi^*(Q) \to 0. 
\end{equation}
with $c_2(\phi^*(Q))=[H]^2 $ and $c_2(\phi^*(S))= 15[H]^2,$ where $[H]$ is the hyper plane 
class. More over it can be easily seen that $ c_1(\phi^*(Q))= 4[H].$ 
Since $E =\phi^*(Q) $ is generated by sections there exists a section
of $E$ which vanishes at exactly one point $p$ of  $ \PP^2$ with multiplicity one.
But then by Cayley-Bacharach theorem(\ref{CB})
the point must be a base point for the complete linear system of the 
line bundle ${\rm det}(E)\otimes K_{\PP^{2}} = \OOO_{\PP^{2}}(1).$
But this is a  contradiction. This proves the required result.
$\hfill{\Box}$
 
 More generally we have the following:
 
 \begin{thm}\label{nonexist} Let $c$ and $a$ be two integers.
  Assume  $c \geq 4 $  and $1\leq a \leq c-2$ or $c^2-c+2 \leq a \leq c^2-1.$
 Then the cohomology class $(a,c^2-a)$  is not represented by a morphism 
 $ \PP^2 \to Gr(2, \CC^{4}).$
 
 \end{thm}
 
 {\it Proof}:   If $\phi :  \PP^2 \to Gr(2, \CC^{4})$ is a morphism such that the cohomology class of the morphism
is $\phi$ is $(a,c^2-a)$ then the cohomology class of the dual morphism (see \ref{re1})
$\phi^{\vee} :  \PP^2 \to Gr(2, \CC^{4})$ is $(c^2-a, a).$ So it is enough to show that there are no morphisms
from $\PP^2$ to $Gr(2, \CC^4)$ such that the cohomolgy class  is given by $(a,c^2-a)$ for
 $1\leq a \leq c-2.$  But note that for any zero dimensional sub scheme $Z$ of $\PP^2$ 
 the natural morphism 
 \[ H^0(\OOO_{\PP^2}(c-3)) \to H^0(\OOO_Z(c-3) \] 
 is surjective, if  length   of $\OOO_Z$ is less than or equal to  $c-2.$ This will imply
 the following: if $Z$ is a zero dimensional sub scheme of $\PP^2$ such that the   
 length  of $\OOO_Z$ is less than or equal to  $c-2,$ then for any $p \in Z$
 \begin{equation}\label{neql}
H^0(\III_Z(c-3)) \neq H^0(\III_{Z-p}(c-3)),
\end{equation} 
where $\III_T$ denotes the ideal sheaf of the sub scheme $T.$
 On the other hand if there exists a morphism  
 $\phi :  \PP^2 \to Gr(2, \CC^{4})$ such that the cohomology class of the morphism
is $\phi$ is $(a,c^2-a)$ then the pull back $E$ of the universal quotient bundle
on $Gr(2, \CC^{4})$ has $c_1(E) =c, c_2(E) = a.$  Since $E$ is
generated by sections we can find section $s$ such that the
scheme $Z$ of zeros of this section is  zero dimensional
and length of $\OOO_Z$ is equal to $a\leq c-2.$ Now by Cayley-Bacharach 
theorem (\ref{CB}) we  must have 
\[ H^0(\III_Z(c-3)) = H^0(\III_{Z-\{p\}}(c-3)) ,\]
for every $p\in Z.$
But this leads to a contradiction to \eqref{neql}. This proves
the reqired result.   $\hfill{\Box}$

 \begin{thm}\label{exist} 
  0) For every $c\geq 1$ the cohomology classes $(0,c^2)$ and $(c^2,0)$
 are represented by  morphisms $\PP^2 \to Gr(2, \CC^{4}).$
  
 1) Let $1\leq c \leq 3$ and $0\leq \ell \leq c^2$ be integers.
 Then  the cohomology class $(\ell,c^2-\ell)$
 is represented by a morphism  $\PP^2 \to Gr(2, \CC^{4}).$
 
 2)  If $c=4$ and $3\leq \ell \leq 13$ then there are morphisms 
$\PP^2 \to Gr(2, \CC^{4})$ representing the cohomology class $(\ell, 16-\ell).$

 3) Let $c \geq 5 $ be an integer. Let $k$ be the largest integer such that 
 $k.c \leq (c^2-3c+2)/2.$ Then for every integer $\ell$ in one of the following
 intervals  $[t(c-3)+2, t.c]$ for $1\leq t \leq k$ or $[(c^2-3c+2)/2+1, c^2/2]$
 there are morphisms $ \PP^2 \to Gr(2, \CC^{4})$ representing the  
 cohomology class $(\ell,c^2-\ell).$ 
 Also, for every such $\ell$ there is a morphism whose cohomology class 
 is $(c^2-\ell, \ell).$
 \end{thm}

{\it Proof:} 
0)  If $f : \PP^2 \to \PP^2$ is a finite morphism then 
$f^*(\OOO_{\PP^2} (1)) = \OOO_{\PP^2} (n)$ for some integer $n>0.$
Then the  degree of the morphism is equal to $n^2.$ Thus we see that
if a cohomology class of the form $(0,a)$ or $(a,0)$ is represented by a
non constant  morphism $f : \PP^2 \to \PP^2$ if and only if $a=n^2$ for a positive integer $n.$

1) If $1\leq c \leq 3$ then for any zero dimensional sub scheme
$Z$ the vector space $H^0(I_Z(c-3))$ is $0.$ Hence by Cayley-Bacharach
theorem(\ref{CB}) there exists a vector bundle $E$ and a section $s$
such that $c_1(E) = c$ and zero scheme of $s$ is $Z.$ It is enough 
to consider the case $\ell \leq c^2/2,$ for if $\phi :  \PP^2 \to Gr(2, \CC^{4})$ such that the 
 cohomology class of the morphism given by $\phi$ is $(\ell,c^2-\ell)$ then 
 the cohomology class of dual morphism (see \ref{re-1}) $\phi^{\vee}$ 
 is $(c^2-\ell, \ell).$
For $\ell\leq c^2/2,$ by considering any reduced sub scheme  $Z$ consisting of $\ell$ points
we get an exact sequence
\[ 0 \to \OOO_{\PP^2} \to E \to I_Z(c) \to 0,\]
where $I_Z$ denotes  the ideal sheaf of $Z,$ where $E$ is the vector bundle 
obtained by Cayley Bacharach theorem. Then we see that $E$ is generated by
sections and hence by $4$ sections. This gives the required morphism $\phi.$ 

2) For $\ell=3$ choose $Z$ a reduced scheme consisting of three points  lying on a line.
Then we see that $Z$ satisfies Cayley - Bacharach conditions with respect to line bundle
$\OOO_{\PP^2}(1).$ Thus there exists vector bundle $E$ with $c_1(det(E))=4[H]$ and 
a section $s$ such that $(s)_0=Z.$ 
Also, We get an exact sequence
\[ 0 \to \OOO_{\PP^2} \to E \to I_Z(4) \to 0,\]
where $I_Z$ denotes  the ideal sheaf of $Z.$  Then we see that $E$ is generated by
sections and hence by $4$ sections. This gives the required morphism $\phi.$ 

For $4 \leq \ell \leq 8$  Let  $Z =\{ P_1, \ldots, P_{\ell} \} $ be a  reduced closed scheme 
consisting of $\ell$ points such that
no three points lie on a line. 
Then we see that $Z$ satisfies Cayley - Bacharach conditions with respect to line bundle
$\OOO_{\PP^2}(1).$ Thus there exists vector bundle $E$ with $c_1(det(E))=4[H]$ and 
a section $s$ such that $(s)_0=Z.$ 
Also, We get an exact sequence
\[ 0 \to \OOO_{\PP^2} \to E \to I_Z(4) \to 0,\]
where $I_Z$ denotes  the ideal sheaf of $Z.$  Then we see that $E$ is generated by
sections and hence by $4$ sections. This gives the required morphism $\phi.$  
For $8\leq \ell \leq 13$ the dual morphism $\phi^{\vee}$  of appropriate $\phi$ above
gives the reqired morphism.
 
3) Let  $t$ and $\ell$ be as in the Theorem. Let  $Z =\{ P_1, \ldots, P_{\ell} \} $ be a  reduced closed scheme 
consisting of $\ell$ points such that
no $r.c+1$ points lie on a curve of degree $r$  for $1\leq r \leq t. $
Note that  $Z$ satisfies Cayley - Bacharach conditions with respect to line bundle
$\OOO_{\PP^2}(c-3).$ Thus there exists vector bundle $E$ with $c_1(det(E))=d[H]$ and 
a section $s$ such that $(s)_0=Z.$ 
Also, We get an exact sequence
\[ 0 \to \OOO_{\PP^2} \to E \to I_Z(d) \to 0,\]
where $I_Z$ denotes  the ideal sheaf of $Z.$  Then we see that $E$ is generated by
sections and hence by $4$ sections. This  morphism $\phi$  corresponds to the 
cohomology class $(\ell, c^2-\ell).$
The dual morphism $\phi^{\vee}$  corresponds to the cohomology class $(c^2-\ell, \ell).$
For  $\ell$  in the interval $[(c^2-3c+2)/2+1, c^2/2]$ let
$Z =\{ P_1, \ldots, P_{\ell} \} $ be a  reduced closed scheme 
consisting of $\ell$ points such that
no $r.c+1$ points lie on a curve of degree $r$  for $1\leq r \leq c-3. $  Then the rest
of the proof is as before. $\hfill{\Box}$

\begin{lem}\label{le22}  Let  $\phi :  \PP^2 \to Gr(2, \CC^{4})$ be a non constant morphism. Assume that 
$\phi^*(Q)$ is decomposable.  
Then the cohomolgy class given by the morphism 
$\phi$ is of the form $(a.b, (a+b)^2-a.b) $  for some non negative integers $a,b$ with $ a+b>0.$
Moreover, for any such tuple  $(a.b, (a+b)^2-a.b),$ there are morphisms 
$\phi :  \PP^2 \to Gr(2, \CC^{4})$ whose corresponding cohomology class is $(a.b, (a+b)^2-a.b).$
\end{lem}

{\it Proof:} By pulling back the universal exact sequence (\eqref{univ2}) by $\phi$
gives the following  exact sequence of vector bundles on $\PP^2:$
\begin{equation}\label{univ7}
 0 \to \phi^*(S) \to \CC^{4} \otimes \OOO_{\PP^{2}} \to \phi^*(Q) \to 0. 
\end{equation}
By assumption $\phi^*(Q) \simeq \OOO_{\PP^{2}}(a)\oplus\OOO_{\PP^{2}}(b).$  Since the bundle 
$\phi^*(Q)$ is quotient of trivial bundle implies we must have $a,b \geq 0.$
On the other hand  $\phi$ is non constant morphism implies $a+b >0.$ Now, it is
easy to see that $c_2(\phi^*(Q)) = a.b$ and $c_2(\phi^*(S)) = (a+b)^2-a.b$ Last assertion of
the Lemma follows from Example(\ref{ex2}). $\hfill{\Box}$

\begin{prop}\label{pr21} Let $X$ and $Y$ be two irreducible projective varieties. Let
$S$ be any irreducible quasi-projective variety and $s_0 \in S$ be a point. Let 
$$F : X \times S \to Y $$
be a morphism. Assume that $F_s :=F|_{X\times s} : X \to Y$ is finite for all $s\in S$
and $F_{s_0}$ is a birational onto its image.  Then there is an open subvariety  $U$ of $S$
such that $s_0 \in U$ and for $s \in U$ the morphism $F_s$ is a birational onto its image.
\end{prop}

{\it Proof:} Consider the morphism $G=F\times Id_S : X \times S \to Y \times S.$
Then the assumption $F_s$ is finite implies the morphism $G$ is finite and proper.
Hence ${\mathcal G}= G_{*}(\OOO_{X\times S})$ is coherent  sheaf of $\OOO_{Y\times S}$ 
modules. Let $Z \subset Y\times S$ be the sub variety on which the sheaf 
$G_{*}(\OOO_{X\times S})$ is supported.  Then clearly the map $p: Z \to S,$ restriction of the natural
projection, is surjective. The section $1 \in \OOO_{X\times S}$ gives an inclusion of 
$ {\OOO_{Z}}$ in ${\mathcal G}.$ Let ${\mathcal F} = {\mathcal G}/{\OOO_{Z}}.$ Let
$Z_1\subset Y\times S $ be the sub variety on which the sheaf ${\mathcal F} $ supported. 
Let $q : Z_1 \to S$ be the natural projection and let
$U = \{ s \in S| {\rm dim}{q^{-1}(s)} < {\rm dim}(X) \} $ then we see that by semi continuity
(see, page 95, Exercise (3.22)  \cite{Ha}),
$U$ is an open subset and is non-empty as $s_0 \in U.$   For $s \in U$ the morphism 
$F_s$ is an isomorphism on $X\times s - G^{-1}(q^{-1}(s).$  Since $G$ is finite 
$G^{-1}(q^{-1}(s)$ is proper closed sub set of $X\times s$ and hence the morphism 
$F_s$  is birational onto its image.
This proves the Proposition.
$\hfill{\Box}$
 
\begin{prop}\label{pr22} Let $a,b$ be two coprime positive  integers. 
Let  
$$f_{\phi} : \PP^2 \to Gr(2, \CC^{4})$$ 
be the morphism given by a surjection 
$$ \phi: \CC^4\otimes\OOO_{\PP^2} \to \OOO_{\PP^2}(a) \otimes  \OOO_{\PP^2}(b).$$
Then for a generic choice of $\phi$ the morphism is birational onto its image.

\end{prop}
 
{\it Proof:} Since the set of surjections $\phi$ is an open subset of 
$$Hom(\CC^4, H^0(\OOO_{\PP^2} (a)\oplus \OOO_{\PP^2} (b)))$$
The result follows from Proposition(\ref{pr21}), if we show the existence of 
one such $f_{\phi}.$  If $\phi$ is given by the matrix
$$
\left[ \begin{array}{cc}
X^a & Z^b\\
Y^a & X^b\\
Z^a & Y^b\\
X^a & 0
\end{array}
\right]
$$
 Since $a,b$ are coprime, it is easy to see that the morphism
$f_{\phi} : \PP^2 \to Gr(2, \CC^{4})$ corresponding to $\phi$
 is birational onto its image.  $\hfill{\Box}$

\begin{rem}\label{re23}  By choosing $a=1$ and $b=n$ and using 
the fact that $Gr(2, \CC^{4})$ is imbedded in $\PP^5$ as a smooth
quadric we get the following:

for every integer $n \geq 1$ 
there are morphisms  
$$f_n: \PP^2  \to  Gr(2, \CC^{4})$$
 such that $f_n$ is one to one onto its image and $f_n^*(\OOO_{Gr(2,\CC^4)}(1)) = 
\OOO_{\PP^2}(n).$ 

\end{rem}

\end{document}